\documentclass{article}
\usepackage{epsf}
\usepackage{latexsym}
\usepackage{amsfonts}

\newtheorem{construction}{\sc Construction}
\newenvironment{proof}{{\sc Proof:}}{~\hfill $\quad\Box$}

\begin{document}

\title{Arbitrary Sectioning of Angles in Taxicab Geometry}

\author{Kevin P. Thompson}

\date{}

\thispagestyle{empty}
\renewcommand\thispagestyle[1]{} 

\maketitle
\begin{abstract}
A construction to arbitrarily section a taxicab angle into an equal number of angles in (pure) taxicab geometry is presented.
\end{abstract}

As humans it is often our nature to investigate new areas by subtly changing something familiar. We like to see the effects of the change and work out the consequences. In geometry, a common approach is to change the distance function, or metric, on a space to obtain a new geometry. One such geometry that can be created in this manner is taxicab geometry.

In Euclidean geometry, we can move from one point to any other point in a straight line and use the usual distance function to calculate how far we traveled. In taxicab geometry, we restrict ourselves to moving along imaginary horizontal and vertical lines like a taxicab driving the streets of a perfectly designed city. This results in a taxicab distance between points $(x_1,y_1)$ and $(x_2,y_2)$ of

\[
d_t=|x_2-x_1|+|y_2-y_1|
\]

\begin{figure}[b]
\centerline{\epsffile{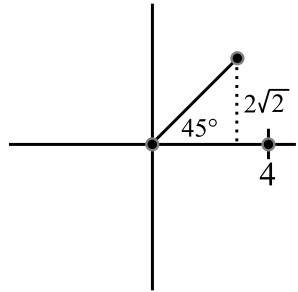}}
\caption{Line segments with different lengths in Euclidean and taxicab geometry.}
\label{distance}
\end{figure}

An immediate result is that line segments with the same Euclidean length often have different taxicab length. The line segments in Figure \ref{distance} have a Euclidean length of 4 but taxicab lengths of 4 and $4\sqrt{2}$. So, the taxicab length of a line segment is dependent on its orientation to the coordinate axes.

To further explore the nature of taxicab geometry, a common approach in introductory non-Euclidean geometry courses is to examine Euclid's axioms and postulates (taxicab geometry fails the SAS congruence postulate) and some of the resulting theorems \cite{Krause, Martin}. A special treatment of Euclid's axioms with respect to taxicab geometry by Dawson \cite{Dawson} in particular examined the angle bisection proposition and the proof of Pierre Wantzel that Euclidean angles cannot be trisected \cite{Wantzel}. It is one aspect Dawson's work that we wish to extend and put a final touch on: a construction to arbitrarily section, or split, a taxicab angle.

\begin{figure}[t]
\centerline{\epsffile{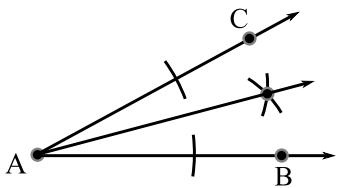}}
\caption{Construction to bisect a Euclidean angle.}
\label{bisection}
\end{figure}

\paragraph{Sectioning a Euclidean Angle}
Using only a straight edge and a compass, the only angle sections that can be created in Euclidean geometry are from repeated bisection. Given an angle $\angle CAB$ (Figure \ref{bisection}), the usual bisection construction proceeds as follows. Use a compass at point $A$ to draw arcs at equal radius intersecting the sides $AB$ and $AC$ of the angle. At the intersection points of these arcs with the sides of the angle, use the compass to draw arcs of equal radius in the interior of the angle. Connecting the intersection of these arcs with the point $A$ will split $\angle CAB$ into two equal angles.

\begin{figure}[b]
\centerline{\epsffile{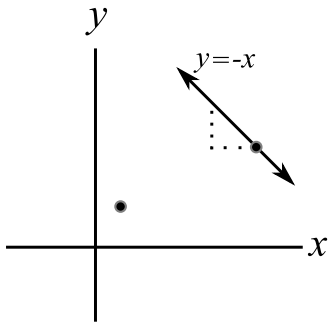}}
\caption{Finding a path in taxicab geometry to keep the distance to a fixed point constant.}
\label{circle_start}
\end{figure}

\paragraph{Taxicab Angles}
While there are multiple ways to choose how an angle is measured in taxicab geometry \cite{Krause, Martin, ThompsonDray}, the most natural is to mimic the Euclidean definition of a radian. For this we need a circle. Having changed the Euclidean metric, even very simple geometric figures take on new forms. The circle needed to define angular measure is one such figure.

\begin{figure}[t]
\centerline{\epsffile{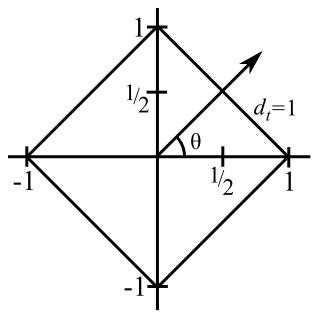}}
\caption{A unit taxicab circle with a 1 t-radian taxicab angle.}
\label{tradian}
\end{figure}

To visualize a circle in taxicab geometry, begin with two points (Figure \ref{circle_start}). If we begin moving from one of the points in a horizontal direction towards the other point, the taxicab distance will get smaller. To keep the distance constant, we must move in the vertical direction to increase the distance and compensate for the horizontal movement. So, distance is kept constant by moving along lines of slope 1 or -1. Extending this observation, we see that the set of all points equidistant from a central point is now a square with its edges oriented $45^{\circ}$ to the horizontal (Figure \ref{tradian}).

Using the unit taxicab circle, we can now define a taxicab radian, or t-radian \cite{ThompsonDray}. A t-radian is an angle that subtends an arc along the unit taxicab circle of length 1 (Figure \ref{tradian}). Note that half a circle covers 4 t-radians. This observation leads to a separate discussion on $\pi$ having the value $\pi_t=4$ in taxicab geometry \cite{Euler}.

So, in sharp contrast to Euclidean geometry, t-radians are measured as length along a straight line and not an arc. Therefore, the procedure for sectioning a taxicab angle reduces to splitting a line segment with slope 1 or -1.

\paragraph{Sectioning a Taxicab Angle}

As noted by Dawson and others, the t-radian can be sectioned using only a compass and a straight edge into an arbitrary number of equal angles. This is a direct result of the t-radian being defined as the length of a line segment and the fact that a line segment can be split into an arbitrary number of equal parts. With taxicab geometry routinely failing to satisfy theorems taken for granted in Euclidean geometry, such a result is a refreshing change of pace!

While Dawson is correct in his statement, and he provides a taxicab construction to bisect a t-radian angle, it was not immediately obvious whether the arbitrary sectioning of an angle could be done with only \emph{taxicab} constructs. Many Euclidean constructions exist to arbitrarily split a line segment, so Euclidean constructs could be used to arbitrarily section a taxicab angle. But, what would this construction look like with only taxicab constructs? Our goal at present is simply to generalize Dawson's construction.

\begin{figure}[t]
\centerline{\epsffile{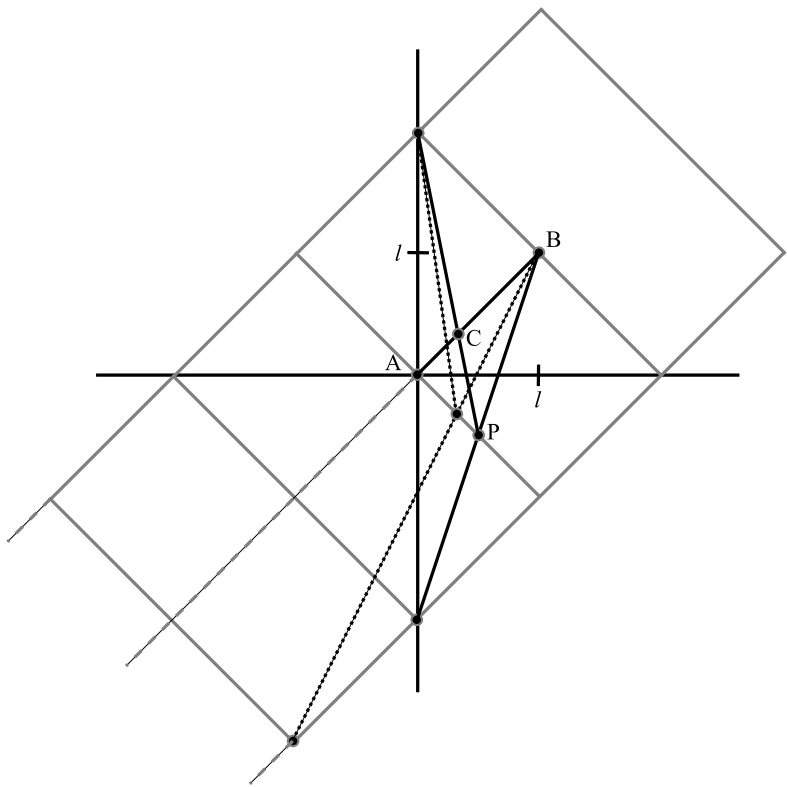}}
\caption{Construction to split a line segment into $n$ equal parts ($n=3$ and $n=4$ shown).}
\label{construction}
\end{figure}

\begin{construction}
To split a line segment with slope 1 into $n\geq3$ equal parts using taxicab constructions
\end{construction}
Given a line segment $\overline{AB}$ of taxicab length $2l$ with slope 1 with $A$ at the origin, construct a taxicab circle of radius $2l$ with center at $B$ (Figure \ref{construction}). Likewise, construct a taxicab circle of radius $2l$ with center at $A$. If $n>3$, at the intersection of the circle around $A$ with the extension of the line segment, construct a circle of radius $2l$. If $n>4$, repeat the last step $n-4$ additional times. Construct a line from the bottom corner of the last circle to point $B$ creating an intersection $P$ with the circle about $B$. Construct a line from $P$ to the top of the circle about $A$ creating an intersection $C$ with the line segment $\overline{AB}$. The taxicab length of line segment $\overline{AC}$ is $2l/n$.

\begin{proof}
Let $A=(0,0)$ and $B=(l,l).$ The bottom corner of the last constructed circle is located at $((3-n)l, (1-n)l)$ so the line from such a point to $B$ is
\[
y=\frac{n}{n-2}(x-l)+l
\]
This line intersects the line $y=-x$ at point $P$ with x-coordinate $x=\frac{l}{(n-1)}$. The line from this intersection to the point $(0,2l)$ at the top of the circle centered at A is
\[
y=(1-2n)x+2l
\]
This line intersects the line $y=x$ at $x=\frac{l}{n}$ giving a taxicab distance from $A$ to this point $C$ of
\[
d(A,C)=\left|\frac{l}{n}-0\right|+\left|\frac{l}{n}-0\right|=\frac{2l}{n}
\]
\end{proof}

\begin{figure}[t]
\centerline{\epsffile{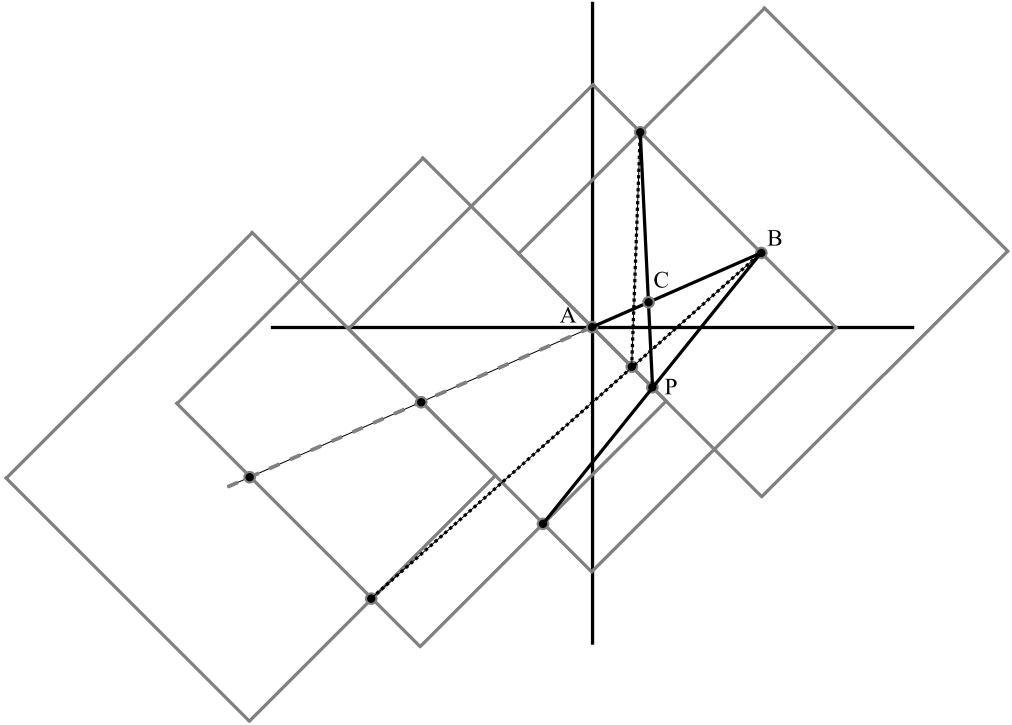}}
\caption{Generalized construction to split a line segment into $n$ equal parts ($n=3$ and $n=4$ shown).}
\label{generalconstruction}
\end{figure}

Dawson's construction for bisection is a special case of the above construction with an adjustment because the bottom corner of the last constructed circle is not defined for $n=2$. This step is merely skipped and the bottom corner of the circle about $A$ is connected to the top of the circle about $B$ which completes the construction.

The presented construction need not be limited to lines of slope 1. This was merely done for the simplicity of calculation and because t-radian angles are defined along diagonal lines. Figure \ref{generalconstruction} illustrates the general construction. Not only have we developed a method to arbitrarily section taxicab angles and line segments, but we have also developed yet another method for arbitrarily splitting Euclidean line segments!

Many times there are multiple constructions that accomplish the same result. The challenge now is to find other taxicab constructions to arbitrarily split a line segment into equal parts!

\end{document}